\documentclass[final]{IEEEtran}

\usepackage[top=1in, left=0.75in, right=0.75in, bottom=1in]{geometry}
\usepackage{amsmath,amssymb,amsthm,graphicx,graphics,times}
\usepackage{algorithm,algorithmic}
\usepackage{color}
\usepackage{ifpdf}
\newtheorem{thm}{Theorem}
\newtheorem{cor}[thm]{Corollary}
\newtheorem{prop}{Proposition}
\newtheorem{lem}[prop]{Lemma}
\newtheorem{mydef}{Definition}

\title{Coherence-Based Performance Guarantees of Orthogonal Matching Pursuit}
\author{
Yuejie Chi$^{\dag}$ and Robert Calderbank$^{\ddag}$ \\
$^{\dag}$Department of Electrical and Computer Engineering, The Ohio State University, USA\\
$^{\ddag}$Department of Electrical and Computer Engineering, Duke University, USA\\
Email: chi@ece.osu.edu, robert.calderbank@duke.edu}

\newcommand{\cI}{{\mathcal I}}
\newcommand{\cG}{{\mathcal G}}
\newcommand{\cH}{{\mathcal H}}

\newcommand{\bbE}{{\mathbb E}}

\newcommand{\snr}{\mathop{\mathsf {SNR}}}

\begin{document}

\maketitle

\thispagestyle{empty}

\begin{abstract}
In this paper, we present coherence-based performance guarantees of Orthogonal Matching Pursuit (OMP) for both support recovery and signal reconstruction of sparse signals when the measurements are corrupted by noise. In particular, two variants of OMP either with known sparsity level or with a stopping rule are analyzed. It is shown that if the measurement matrix $X\in\mathbb{C}^{n\times p}$ satisfies the strong coherence property, then with $n\gtrsim\mathcal{O}(k\log p)$, OMP will recover a $k$-sparse signal with high probability. In particular, the performance guarantees obtained here separate the properties required of the measurement matrix from the properties required of the signal, which depends critically on the minimum signal to noise ratio rather than the power profiles of the signal. We also provide performance guarantees for partial support recovery. Comparisons are given with other performance guarantees for OMP using worst-case analysis and the sorted one step thresholding algorithm.

\end{abstract}


\section{Introduction}

Sparse signal processing is a fundamental task in many applications involving high-dimensional data. Many of the recent advances in Compressive Sensing (CS) \cite{Donoho,CandesTao} have been centered around reconstructing a sparse signal $\beta\in\mathbb{C}^p$ with a few nonzero entries, from a number of linear measurements that is much smaller than the signal dimension, possibly corrupted by noise $\eta$, given as $y=X\beta+\eta$, where $X\in\mathbb{C}^{n\times p}$ is the measurement matrix. Two main classes of algorithms have been successfully applied, one is convex optimization based algorithms such as Basis Pursuit (BP) \cite{Chen2001}, the other one is greedy pursuit based algorithms such as Orthogonal Matching Pursuit (OMP) \cite{Tropp_OMP}. In particular, the latter class is appealing and competitive in practice due to its simplicity and low computational cost \cite{Kunis_Rauhut_2007}.

The performance of OMP can be characterized either in a worst-case sense or in an average (probabilistic) sense. Define the worst coherence of the measurement matrix $X$ as $\mu =\max_{i\neq j} |\langle x_i,x_j \rangle |$, where $x_i$ denotes the $i$th unit-norm column of $X$. It is shown in \cite{Tropp_OMP} that if $\mu<\frac{1}{2k-1}$, then OMP recovers any $k$-sparse vector $\beta$ from the noiseless measurement $y=X\beta$, and this result is confirmed to be sharp in \cite{Cai2010a}. It is further studied in \cite{CaiWang} that given the amplitudes of the nonzero entries of $\beta$ are not too small, OMP recovers the support of the signal from noisy observations. From the Welch bound \cite{Welch} which gives $\mu\gtrsim\mathcal{O}(n^{-1/2})$, in order to recover \textit{all} $k$-sparse vectors, the number of measurements is required to satisfy $n\gtrsim\mathcal{O}(k^2)$. It is demonstrated in \cite{OMP-RIP} that, if the Restricted Isometry Property (RIP) of $X$ of order $k+1$ is satisfied with constant smaller than $1/(3\sqrt{k})$, then OMP recovers any $k$-sparse signal, again from the noiseless measurements $y=X\beta$. Given a random matrix satisfies the RIP of order $k$ with constant $\delta_k$ provided that $n\gtrsim\mathcal{O}(k\log(p/k)/\delta_k^2)$, then the number of measurements is required $n\gtrsim\mathcal{O}(k^2\log (p/k))$. These results all suffer from what is called ``square-root bottleneck''. Alternatively, instead of aiming to recover all $k$ sparse signals using OMP, it is shown that for a fixed sparse vector $\beta$, a randomly drawn measurement matrix $X$ from i.i.d. normal entries can recover $\beta$ with high probability with $n\gtrsim\mathcal{O}(k\log p)$ measurements from the noiseless measurements $y=X\beta$. However, it does not provide a practical way to design or validate the usefulness of a measurement matrix.

In this paper, we aim to use the strong coherence properties proposed in \cite{BCJ2010} to capture the performance of OMP. The strong coherence properties require that the worst-case coherence of the measurement matrix $X$ is sufficiently small, and also that the average coherence, defined as $\nu =\frac{1}{p-1} \max_i  | \sum_{j:j\neq i}\langle x_i,x_j \rangle |$, is small relative to the worst-case coherence. Different from the worst-case sense, we aim to provide conditions on the measurement matrix $X$ that are easily verifiable in contrast to RIP, and that succeeds in support recovery and reconstruction of sparse signals with high probability in the presence of noise. We show that if $X$ satisfied the strong coherence property, then with $n\gtrsim\mathcal{O}(k\log p)$ measurements, OMP recovers a $k$-sparse signal with high probability. In particular, the performance of OMP depends on the smallest signal-to-noise ratio $\snr_{\min}$ determined by the smallest nonzero entry instead of the relative strengths of the nonzero entries of the signal. 

The rest of the paper is organized as follows. Section~\ref{backgrounds} introduces two coherence parameters of the measurement matrix and provides a detailed description of OMP. Section~\ref{main} gives the main theorems on the performance guarantees of OMP for support recovery and signal reconstruction under the strong coherence property. Section~\ref{prep} prepares for the proof and Section~\ref{proof} proves the main theorems. Finally we conclude with discussions in Section~\ref{conclusion}.

\section{Two Fundamental Coherence Parameters} \label{backgrounds}

Suppose we are given a measurement vector as $y= X\beta+\eta$, where $X$ is an $n\times p$ unit-norm measurement matrix, $\beta\in\mathbb{C}^p$ is a $k$-sparse vector, $y\in\mathbb{C}^n$ is the measurement and $\eta\in\mathbb{C}^n$ is the noise. We define two coherence properties of $X=[x_1,\cdots,x_p]$ below. The first is worst-case coherence:
\begin{equation*}
\mu =\max_{i\neq j} |\langle x_i,x_j \rangle |=\| X^HX-I \|_\infty,
\end{equation*}
which captures the correlation between different columns of $X$. The second is average coherence: 
\begin{align*}
\nu & =\frac{1}{p-1} \max_i \Big | \sum_{j:j\neq i}\langle x_i,x_j \rangle \Big |=\frac{1}{p-1} \|(X^HX-I)1\|_\infty.
\end{align*}
which captures the average correlation between one column of $X$ and the remaining columns of $X$.

We say a measurement matrix $X$ satisfies the strong coherence property if the equation below holds:
\begin{equation} \label{CP}
\mu \leq \frac{1}{240\log p}, \quad\quad \nu \leq \frac{\mu}{\sqrt{n}}.
\end{equation}

It is known that Gaussian random matrices satisfy the strong coherence property with high probability as long as $n\gtrsim\mathcal{O}((\log p)^4)$ \cite{BCD2011}. Several families of deterministic matrices are also known to satisfy the strong coherence property, including Gabor frames \cite{BCJ2010}, Kerdock code sets \cite{CalderbankLasso}, and Delsarte-Goethals code sets \cite{CalderbankLasso}. 
%show some nuber
Notice that the condition on average coherence $ \nu \leq \frac{\mu}{\sqrt{n}}$ can be achieved with essentially no cost via``wiggling'', i.e. flipping the signs of the columns of $X$ \cite{BCD2011}. The ``wiggling'' procedure doesn't change $\mu$ and $\|X\|_2$.

The goal of this paper is to present performance guarantees on OMP for both support recovery and signal reconstruction under the assumption that $X$ satisfies the strong coherence property, when the measurements are corrupted by noise. In particular, two variants of OMP that differ in the way they terminate the iterations are analyzed, i.e. Algorithm~\ref{OMP-k} with known sparsity level, and Algorithm~\ref{OMP-stop} with a stopping rule.

\begin{algorithm}
\caption{OMP with a fixed number of iterations}
\label{OMP-k}
\begin{algorithmic}[1]
\STATE Input: an $n\times p$ matrix $X$, a vector $y\in\mathbb{C}^n$, and a sparsity level $k$
\STATE Output: an estimate $\hat{\mathcal{S}}$ of the true model $\mathcal{S}$
\STATE Initialization: $\hat{\mathcal{S}}_0:=$ empty set, residual $r_0=y$
\FOR{$t:=1:k$}
\STATE $f:=X^Hr_{t-1}$\label{correlation}
\STATE $i:=\arg\max_{j} |f_j|$ 
\STATE $\hat{\mathcal{S}}_t:=\hat{\mathcal{S}}_{t-1}\cup\{j\}$
\STATE $r_t := y - X_{\hat{\mathcal{S}}_t}(X_{\hat{\mathcal{S}}_t}^HX_{\hat{\mathcal{S}}_t})^{-1} X_{\hat{\mathcal{S}}_t}^Hy$
\ENDFOR
\STATE $\hat{\mathcal{S}} :=\hat{\mathcal{S}}_k$
\end{algorithmic} 
\end{algorithm} 

\begin{algorithm}
\caption{OMP with a stopping rule}
\label{OMP-stop}
\begin{algorithmic}[1]
\STATE Input: an $n\times p$ matrix $X$, a vector $y\in\mathbb{C}^n$, and a threshold level $\delta$
\STATE Output: an estimate $\hat{\mathcal{S}}$ of the true model $\mathcal{S}$
\STATE Initialization: $\hat{\mathcal{S}}_0:=$ empty set, residual $r_0=y$, set the iteration counter $t=1$
\WHILE{$\|X^Hr_{t-1} \|_\infty >  \delta $}
\STATE $f:=X^Hr_{t-1}$
\STATE $i:=\arg\max_{j} |f_j|$
\STATE $\hat{\mathcal{S}}_t:=\hat{\mathcal{S}}_{t-1}\cup\{j\}$
\STATE $r_t := y - X_{\hat{\mathcal{S}}_t}(X_{\hat{\mathcal{S}}_t}^HX_{\hat{\mathcal{S}}_t})^{-1} X_{\hat{\mathcal{S}}_t}^Hy$
\STATE $t:=t+1$
\ENDWHILE
\STATE $\hat{\mathcal{S}} :=\hat{\mathcal{S}}_{t-1}$
\end{algorithmic} 
\end{algorithm}

\section{Performance Guarantees of OMP} \label{main}

We define the minimum-to-average ratio $\mathsf{MAR}$ and the $t$th-largest-to-average ratio $\mathsf{LAR}_{(t)}$ of the signal $\beta$ respectively as
\begin{align}
\mathsf{MAR} & = \frac{|\beta|_{\min}^2}{\|\beta\|_2^2/k}, \quad \mathsf{LAR}_{(t)} = \frac{|\beta|_{(t)}^2}{\|\beta\|_2^2/k},
\end{align}
where $|\beta|_{(t)}$ is the $t$th largest absolute value of $\beta$, $|\beta|_{\min}$ is the smallest nonzero absolute value of $\beta$, and $k$ is the sparsity level of $k$. The signal-to-noise ratio $\snr$ and minimum signal-to-noise ratio $\snr_{\min}$ are defined respectively as
\begin{equation}
\snr  = \frac{\|\beta\|_2^2}{\mathbb{E}\|\eta\|_2^2}, \quad \mathsf{SNR}_{\min}  =  \frac{|\beta|_{\min}^2}{\mathbb{E}\|\eta\|_2^2/k}. 
\end{equation}

\subsection{Performance Guarantee for Support Recovery}
We have the following theorem for OMP with knowledge of sparsity level $k$ in Algorithm~\ref{OMP-k}.
\begin{thm} \label{thm-omp-k}
Suppose $X$ satisfies the strong coherence property for any $p\geq 128$, and $\eta\sim\mathcal{CN}(0,\sigma^2I)$. If the sparsity level of $\beta$ satisfies 
\begin{equation} \label{k-cond}
k\leq \min\left\{ \frac{p}{c_2^2\|X\|_2^2\log p}, \frac{1}{c_1^2\mu^2\log p} \right\}
\end{equation}
for $c_1=50\sqrt{2}$ and $c_2=104\sqrt{2}$, and its nonzero entries satisfy
\begin{equation} \label{k-condition}
|\beta|_{(t+1)} >\frac{2\sigma\sqrt{(1+\alpha)\log p}}{1-c_1\mu\sqrt{(k-t) \log p} },  
 \end{equation}
or write differently, that is 
\begin{equation} \label{lar-condition}
 \mathsf{LAR}_{(t+1)} >\frac{4(1+\alpha)}{(1-c_1\mu\sqrt{(k-t) \log p})^2 } \cdot \left( \frac{k \log p}{n\snr } \right),  
 \end{equation}
for $0\leq t\leq k-1$ and $\alpha\geq1$, then the OMP algorithm in Algorithm~\ref{OMP-k} successfully finds the support of $\beta$ with probability at least $1-k(p^{\alpha}\pi)^{-1}-2p^{-2\log 2}-4p^{-1}$. 
\end{thm}

For the OMP algorithm with a stopping rule in Algorithm~\ref{OMP-stop}, we have the following theorem.
\begin{thm} \label{thm-omp-stop}
If $X$ satisfies the strong coherence property for any $p\geq 128$, and $\eta\sim\mathcal{CN}(0,\sigma^2I)$. If the sparsity level of $\beta$ satisfies \eqref{k-cond} and its nonzero entries satisfy \eqref{lar-condition} for $\alpha\geq 1$, and choose $\delta=\sigma\sqrt{(1+\alpha)\log p}$,  then the OMP algorithm in Algorithm~\ref{OMP-stop} successfully finds the support of $\beta$ with probability at least $1-(k+1)(p^{\alpha}\pi)^{-1}-2p^{-2\log 2}-4p^{-1}$ in $k$ iterations. 
\end{thm}

%Remark: the decay law in the signal amplitude required in Theorem~\ref{thm-omp-k} and Theorem~\ref{thm-omp-stop} is at least
% \mathsf{LAR}_{(t+1)} 

Since $\mathsf{MAR}\leq \mathsf{LAR}_{(t+1)}$ for all $0\leq t\leq k-1$, we have the following corollary respectively from Theorem~\ref{thm-omp-k} and Theorem~\ref{thm-omp-stop}.
\begin{cor}  \label{cor-k-mar}
If $X$ satisfies the strong coherence property for any $p\geq 128$, and $\eta\sim\mathcal{CN}(0,\sigma^2I)$. If the sparsity level of $\beta$ satisfies \eqref{k-cond}, and it satisfies
\begin{equation} \label{mar-condition}
 \mathsf{MAR} >\frac{4(1+\alpha)}{(1-c_1\mu\sqrt{k\log p})^2 } \cdot \left( \frac{k \log p}{n\snr } \right), 
 \end{equation}
for $\alpha\geq 1$, then the support of $\beta$ is successfully recovered by OMP with probability at least $1-k(p^{\alpha}\pi)^{-1}-2p^{-2\log 2}-4p^{-1}$ using Algorithm~\ref{OMP-k}, and with probability at least $1-(k+1)(p^{\alpha}\pi)^{-1}-2p^{-2\log 2}-4p^{-1}$ by Algorithm~\ref{OMP-stop} with $\delta=\sigma\sqrt{(1+\alpha)\log p}$ in $k$ iterations. 
\end{cor}

Let $\theta=c_1\mu\sqrt{k\log p}\in(0,1)$, then \eqref{mar-condition} implies  the sparsity level $k$ satisfies
\begin{equation} \label{k-extra-cond}
k<\frac{(1-\theta)^2}{4(1+\alpha)} \cdot \frac{\snr_{\min}}{\log p}.
\end{equation}
Combining with \eqref{k-cond}, we have
\begin{align} \label{k-all-cond}
k <\max_{0<\theta<1} \min\Big\{ &\frac{n(1-\theta)^2\snr_{\min}}{4(1+\alpha)\log p}, \frac{\theta^2}{c_1^2\mu^2\log p}, \nonumber \\ 
&\quad \frac{p}{c_2^2\|X\|_2^2\log p} \Big\},
\end{align}
where the first term is determined by $\snr_{\min}$, which is signal dependent; the second term and the third term are determined by the worst-case coherence and the spectral norm of the measurement matrix $X$. If $X$ is a tight frame, $\|X\|_2^2=p/n$, the third term becomes $k< \mathcal{O}(n/\log p)$. From the Welch bound $\mu$ is lower bounded as $\mu\gtrsim\mathcal{O}(n^{-1/2})$, we write the worst-case coherence as $\mu=c_3n^{-1/\gamma}$ for some $c_3>0$ and $\gamma\geq 2$. Therefore the maximum sparsity level is determined by the second term in \eqref{k-all-cond}, yielding $k\lesssim \mathcal{O}((n/\log p)^{2/\gamma})$, and when $\gamma=2$ this gives $k\lesssim \mathcal{O}(n/\log p)$.  In particular, the sparsity level $k$ doesn't depend on the profile of signal strength of $\beta$, i.e. $\mathsf{MAR}$ of the signal.

We have another corollary on partial recovery.
\begin{cor} \label{partial}
If $X$ satisfies the strong coherence property for any $p\geq 128$, and $\eta\sim\mathcal{CN}(0,\sigma^2I)$. If the sparsity level of $\beta$ satisfies \eqref{k-cond}, and its nonzero entries satisfy \eqref{lar-condition} for $0\leq t\leq k'-1\leq k-1$ and $\alpha\geq 1$, then the OMP algorithm in both Algorithm~\ref{OMP-k} and Algorithm~\ref{OMP-stop} successfully selects $k'$ entries from the support of $\beta$ with probability at least $1-k'(p^{\alpha}\pi)^{-1}-2p^{-2\log 2}-4p^{-1}$. 
\end{cor}

It is worth noting that it is not necessarily the support of the $k'$-largest entries that is recovered from the first $k'$ iterations. The next corollary provides the condition on detecting the $k'$-largest entries.

\begin{cor} \label{partial-largest}
If $X$ satisfies the strong coherence property for any $p\geq 128$, and $\eta\sim\mathcal{CN}(0,\sigma^2I)$. If the sparsity level of $\beta$ satisfies \eqref{k-cond}, and its nonzero entries satisfy
\begin{equation} \label{decay} 
|\beta|_{(t+1)} > \frac{|\beta|_{(t+2)}+ 2\sigma\sqrt{(1+\alpha)\log p}}{1-c_1\mu\sqrt{(k-t) \log p} }, 
\end{equation}
for $0\leq t\leq k'-1\leq k-1$ and $\alpha\geq 1$, then the OMP algorithm in both Algorithm~\ref{OMP-k} and Algorithm~\ref{OMP-stop} successfully selects $k'$ largest entries from the support of $\beta$ with probability at least $1-k'(p^{\alpha}\pi)^{-1}-2p^{-2\log 2}-4p^{-1}$. 
\end{cor}

\subsection{Performance Guarantees for Signal Reconstruction}
Furthermore, we could reconstruct the amplitude of the signal $\beta$ by first reconstructing the amplitude on the detected support $\Pi$ via least-squares estimation as 
$$\hat{z}=X_{\Pi}^{\dag}y,$$
then $\hat{\beta}$ is obtained by filling in the zero entries. We have the following theorem.

\begin{thm} \label{reconstruction}
If $X$ satisfies the strong coherence property for any $p\geq 128$, and $\eta\sim\mathcal{N}(0,\sigma^2I)$. If the sparsity level of $\beta$ satisfies \eqref{k-cond} and its nonzero entries satisfy \eqref{lar-condition} for $0\leq t\leq k-1$ and $\alpha\geq 1$, then the $\ell_2$ norm error of the signal reconstructed by least-squares estimation on the support recovered by the OMP algorithm satisfies
$$ \| \hat{\beta}-\beta\|_2^2 \leq 4(1+\alpha)k\sigma^2\log p $$
with probability at least $1-k(p^{\alpha}\pi)^{-1}-2p^{-2\log 2}-4p^{-1}$ using Algorithm~\ref{OMP-k}, and with probability at least $1-(k+1)(p^{\alpha}\pi)^{-1}-2p^{-2\log 2}-4p^{-1}$ using Algorithm~\ref{OMP-stop}.
\end{thm}

\subsection{Comparison with Other Results}
We now compare our bound with the performance guarantee of OMP for support recovery provided in \cite{CaiWang}, which we have modified slightly for complex Gaussian noise. In order to select exactly the correct support with probability at least $1-(k+1)(p^{\alpha}\pi)^{-1}$ for the OMP Algorithm~\ref{OMP-stop} with the stopping rule $\delta=\sigma\sqrt{(1+\alpha)\log p}$, the signal $\beta$ needs to satisfy
\begin{equation}
 \mathsf{MAR} >\frac{4(1+\alpha)}{(1-(2k-1)\mu)^2 } \cdot \left( \frac{k \log p}{n\snr } \right),
\end{equation}
therefore the sparsity level of $\beta$ satisfies
\begin{equation*} \label{k-cond-1}
k <\max_{0<\theta<1} \min \left\{ \frac{n(1-\theta)^2\snr_{\min}}{4(1+\alpha)\log p}, \frac{1}{2}+\frac{\theta}{2\mu} \right\}.
\end{equation*}
The first term is the same as that in \eqref{k-all-cond}, but the second term gives $k\sim\mathcal{O}(\mu^{-1})$, therefore $k\lesssim\mathcal{O}(n^{1/\gamma})$. We achieved a much tighter bound \eqref{k-all-cond} by sacrificing the probability of success to $1-(k+1)(p^{\alpha}\pi)^{-1}-2p^{-2\log 2}-4p^{-1}$.

We also compare with the performance guarantee of the Sorted One Step Thresholding (SOST) algorithm analyzed in \cite{BCJ2010}, which outputs the index set of the $k$-largest entries in absolute values of  $f=X^Hy$ from line \ref{correlation} in Algorithm~\ref{OMP-k}. By rephrasing Theorem 4 in \cite{BCJ2010}, in order to select the correct support with probability at least $1-6p^{-1}$, the sparsity level of $\beta$ satisfies
\begin{align} \label{k-sost}
k <\max_{0<\theta<1} \min\Big\{ &\frac{n(1-\theta)^2\snr_{\min}}{16\log p}, \frac{\theta^2}{800\mu^2\log p}\cdot\frac{1}{\mathsf{MAR}}, \nonumber \\ 
&\quad \frac{n}{2\log p} \Big\}.
\end{align}
Compare to \eqref{k-all-cond}, we see that when $X$ is not a tight frame, it is the third term that degrades the performance guarantees. On the other hand, the SOST algorithm performs poorly when the $\mathsf{MAR}$ is of the signal is much smaller than $1$, as seen from the second term in \eqref{k-sost}.

\section{Preparations for Proof} \label{prep}
\subsection{Statistical Orthogonality Condition (StOC)}
The Statistical Orthogonality Condition (StOC) for a measurement matrix $X$ is first introduced in \cite{BCJ2010}.

\begin{mydef} Let $\bar{\Pi}=(\pi_1,\ldots,\pi_p)$ be a  random permutation of $\{1,\ldots,p\}$, and define $\Pi=(\pi_1,\ldots,\pi_k)$ and $\Pi^c=(\pi_{k+1},\ldots,\pi_{p})$ for any $k\leq p$. Then the matrix $X$ is said to satisfy the $(k,\epsilon,\delta)$-StOC, if there exist $\epsilon$, $\delta\in[0,1)$ such that the inequalities 
\begin{align} 
\| (X_\Pi^HX_\Pi-I)z \|_{\infty} & \leq\epsilon \|z\|_2, \label{StOC-1}\\
\| X_{\Pi^c}^HX_\Pi z \|_{\infty} & \leq\epsilon \|z\|_2, \label{StOC-2}
\end{align}
hold for every fixed $z\in\mathbb{C}^k$ with probability exceeding $1-\delta$, with respect to $\bar{\Pi}$. 
\end{mydef}

We have the following proposition rephrased from \cite{BCJ2010} stating that the StOC is satisfied with high probability if $X$ satisfies the strong coherence property. 
\begin{prop}[\cite{BCJ2010}] \label{StOC}
If the measurement matrix $X$ satisfies the strong coherence property, then it satisfies $(k,\epsilon,\delta)$-StOC for $k\leq n/(2\log p)$, with $\epsilon=10\mu\sqrt{2\log p}$ and $\delta\leq 4p^{-1}$.
%any $\epsilon\in[0,1)$ with $\delta\leq 4p\exp\left(-\frac{(\epsilon-\sqrt{k}\nu)^2}{16(2+a^{-1})^2\mu^2}\right)$ for any $a\geq 1$ as long as $k\leq \min\{\epsilon^2\nu^{-2}, (1+a)^{-1}p\}$.
\end{prop}

If \eqref{StOC-1} and \eqref{StOC-2} hold for a realization of permutation $\bar{\Pi}$, then for $t\leq k$, let $\Pi_t=(\pi_1,\ldots,\pi_{t})$ and $\Pi_t^c=(\pi_{t+1},\ldots,\pi_{k})$, if \eqref{StOC-1} and \eqref{StOC-2} hold for every $z\in\mathbb{C}^k$, so that $\Pi_t\cup\Pi_t^c=\Pi$ and $\Pi_t\cap\Pi_t^c=\emptyset$. For every $z \in\mathbb{C}^{t}$, we have
\begin{align*}
\left \| \begin{bmatrix}
 X_{\Pi_t}^HX_{\Pi_t}-I_t  & X_{\Pi_t}^H X_{\Pi_t^c} \\
  X_{\Pi_t^c}^H X_{\Pi_t} &  X_{\Pi_t^c}^HX_{\Pi_t^c}-I_{k-t}
  \end{bmatrix} \begin{bmatrix}
  z \\
  0
  \end{bmatrix} \right\|_{\infty} & \leq\epsilon \|z\|_2, 
\end{align*}
from \eqref{StOC-1}, therefore
\begin{align*}
\| (X_{\Pi_t}^HX_{\Pi_t}-I_t)z \|_{\infty} & \leq\epsilon \|z\|_2, \\
 \| X_{\Pi_t^c}^HX_{\Pi_t} z \|_{\infty} & \leq\epsilon \|z\|_2. 
\end{align*}
Moreover, from \eqref{StOC-2} we have
\begin{align*}
\| X_{\Pi^c}^HX_{\Pi_t} z \|_{\infty} & = \left\| \begin{bmatrix}
  X_{\Pi^c}^H X_{\Pi_t} &  X_{\Pi^c}^HX_{\Pi_t^c}
  \end{bmatrix} \begin{bmatrix}
  z \\
  0
  \end{bmatrix} \right\|_{\infty} \leq\epsilon \|z\|_2.
\end{align*}

%If we let $k\leq n/(2\log p)$ and fix $\epsilon=10\mu\sqrt{2\log p}$, then $X$ satisfy the $(k,\epsilon,\delta)$-StOC with $\delta\leq 4p^{-1}$ provided that $X$ satisfies the coherence property.

\subsection{Conditioning of random submatrices}
We need the following proposition that shows a random submatrix of $X$ is well-conditioned with high probability, which is essentially due to Tropp \cite{Tropp}, and first presented in the form below by Cand\`es and Plan \cite{CandesPlan}. 

\begin{prop}[\cite{Tropp,CandesPlan}] \label{submatrix}
Let $\bar{\Pi}=(\pi_1,\ldots,\pi_p)$ be a  random permutation of $\{1,\ldots,p\}$, and define $\Pi=(\pi_1,\ldots,\pi_k)$ for any $k\leq p$. Then for $q=2\log p$ and $k\leq p/(4\|X\|_2^2)$, we have
\begin{align}
 \quad  & \left(\bbE\left[ \|X_{\Pi}^HX_{\Pi} - I \|_2^q \right] \right)^{1/q} \nonumber \\
& \leq 2^{1/q} \left(30\mu\log p+13\sqrt{\frac{2k\|X\|_2^2\log p}{p}} \right).
\end{align}
with respect to the random permutation $\bar{\Pi}$. 
%Moreover, for the same value of $q$,
%\begin{equation}
%\left(\bbE\left[ \|X_{\Pi}^HX_{\Pi^c} \|_\infty^q \right] \right)^{1/q} \leq 2^{1/q} \left(4\mu\sqrt{\log p}+13\sqrt{\frac{k\|X\|_2^2}{p}} \right).
%\end{equation}
\end{prop}

The following proposition \cite{CandesPlan} states a probabilistic bound on the extreme singular values of a random submatrix of $X$, by applying Markov's inequality 
$$ \Pr\left( \|X_{\Pi}^HX_{\Pi}-I \|_2\geq 1/2\right) \leq 2^q \bbE\left[ \|X_{\Pi}^HX_{\Pi} - I \|_2^q\right] $$
to Proposition~\ref{submatrix}.
\begin{prop}[\cite{CandesPlan}] \label{submat}
Let $\bar{\Pi}=(\pi_1,\ldots,\pi_p)$ be a  random permutation of $\{1,\ldots,p\}$, and define $\Pi=(\pi_1,\ldots,\pi_k)$ for any $k\leq p$. Suppose that $\mu(X)\leq 1/(240\log p)$ and $k\leq p/(c_2^2\|X\|_2^2\log p)$ for numerical constant $c_2=104\sqrt{2}$, then we have
$$ \Pr\left( \|X_{\Pi}^HX_{\Pi}-I \|_2\geq 1/2\right) \leq 2p^{-2\log 2}.$$
\end{prop}
Notice that $\|X_{\Pi}^HX_{\Pi}-I \|_2= \max \{\lambda_{\max}(X_{\Pi}^HX_{\Pi})-1, 1-\lambda_{\min}(X_{\Pi}^HX_{\Pi})\}$, where $\lambda_{\max}(X_{\Pi}^HX_{\Pi})$ and  $\lambda_{\min}(X_{\Pi}^HX_{\Pi})$ are the maximum and minimum eigenvalues of $X_{\Pi}^HX_{\Pi}$, i.e. all the eigenvalues of $X_{\Pi}^HX_{\Pi}$ are bounded in $[1/2, 3/2]$. If for a realization of permutation $\bar{\Pi}$, $ \|X_{\Pi}^HX_{\Pi}-I \|_2\geq 1/2$, we have
$$\| (X_{\Pi}^HX_{\Pi})^{-1}\|_2\leq 2~~\mbox{and}~~\|X_{\Pi}(X_{\Pi}^HX_{\Pi})^{-1}\|_2\leq\sqrt{2}.$$
% hold under $\cG_2$.
Moreover, for $t\leq k$ and $\Pi_t=(\pi_1,\ldots,\pi_{t})$, we have 
$$ \|X_{\Pi_t}^HX_{\Pi_t}-I_t \|_2 \geq 1/2, $$
since the eigenvalues of $X_{\Pi_t}^HX_{\Pi_t}$ are majorized by the eigenvalues of $X_{\Pi}^HX_{\Pi}$.

\subsection{Correlated Gaussian Noise}

%Let the noise $\eta$ be a random vector with i.i.d. $\mathcal{N}(0, \sigma^2)$ entries. 
%We have $\hat{\eta}=X^H\eta$ satisfies $\| X^H\eta \|_\infty \leq 2\sqrt{\sigma^2 \log p}$ with probability $1-2(p\sqrt{2\pi \log p} )^{-1}$. Define the corresponding event 
%$\cG_3=\{\| X^H\eta \|_\infty \leq 2\sqrt{\sigma^2 \log p}\}$.
Let $P\in\mathbb{C}^{n\times n}$ be a projection matrix such that $P^2=P$. Since $\eta\sim\mathcal{CN}(0,\sigma^2I_n)$ is i.i.d. complex Gaussian noise, $X^HP\eta\sim\mathcal{CN}(0,\sigma^2X^HPX)$ is also Gaussian distributed, but is correlated with covariance matrix $\sigma^2X^HPX$. We want to bound $\Pr(\| X^HP\eta \|_\infty\geq \tau)$ for some $\tau>0$. First, we need the Sidak's lemma \cite{Sidak} below.
\begin{lem}[Sidak's lemma]
Let $[X_1,\cdots, X_n]$ be a vector of random multivariate normal variables with zero means, arbitrary variances $\sigma_1^2$, $\cdots$, $\sigma_n^2$ and
and an arbitrary correlation matrix. Then, for any positive numbers $c_1, \cdots, c_n$, we have
$$ \Pr( |X_1| \leq c_1, \cdots,  |X_n| \leq c_n) \geq \prod_{i=1}^n  \Pr( |X_i| \leq c_i). $$
\end{lem}
Since $X^HP\eta\sim\mathcal{CN}(0,\sigma^2X^HPX)$, then each $x_i^HP\eta\sim\mathcal{CN}(0,\sigma_i^2)$, where $\sigma_i^2=\sigma^2 x_i^HPx_i\leq \sigma^2$. Then
$$ \Pr(| x_i^HP\eta | \leq \tau) =1-\frac{1}{\pi}e^{-\tau^2/\sigma_i^2}\geq 1-\frac{1}{\pi}e^{-\tau^2/\sigma^2}. $$
Following Sidak's lemma, for $\tau>0$ we have
\begin{align*}
\Pr(\| X^HP\eta \|_\infty \leq \tau)& \geq \prod_{i=1}^p \Pr(| x_i^HP\eta | \leq \tau) \\
& \geq (1-\frac{1}{\pi}e^{-\tau^2/\sigma^2})^p \\
& \geq 1-\frac{p}{\pi}e^{-\tau^2/\sigma^2},
\end{align*}
provided the RHS is greater than zero. We have the proposition below.

%Since $|x_i^HPx_i|\leq \|x_i\|_2^2= 1$, we could generalize Proposition 4 in \cite{TroppGilbert} for complex Gaussian variables as below.
\begin{prop} \label{noise}
Let $\eta$ be a random vector with i.i.d. $\mathcal{CN}(0, \sigma^2)$ entries, $P$ be a projection matrix, and $X$ be a unit-column matrix, then for $\tau>0$ we have
\begin{align*}
\Pr(\| X^HP\eta \|_\infty \leq \tau)& \geq 1-\frac{p}{\pi}e^{-\tau^2/\sigma^2},
\end{align*}
provided the RHS is greater than zero.
\end{prop}

Now let $\tau=\sigma\sqrt{ (1+\alpha)\log p}$ for $\alpha\geq 1$, we have
$$ \Pr \{ \|X^HP\eta \|_\infty \leq \sqrt{\sigma^2(1+\alpha)\log p} \} \geq 1-(p^{\alpha}\pi)^{-1}. $$

\section{Proof of Main Results} \label{proof}

We first write the data vector $\beta$ as $\beta=P_{\Pi}z$, where $z\in\mathbb{C}^{k}$ is a deterministic vector, and $P_{\Pi}\in\mathbb{R}^{p\times k}$ is a partial identity matrix composed of columns indexed by $\Pi$. Then the measurement vector can be written as
$$y=X\beta+\eta = XP_{\Pi}z+\eta=X_{\Pi}z+\eta, $$
where $X_{\Pi}$ denotes the submatrix of $X$ composed of columns indexed by $\Pi$.

We note that in OMP, the residual $r_t$, $t=0,\cdots, k-1$ is orthogonal to the selected columns in previous iterations, so in each iteration a new column will be selected. Define a subset $\Pi_t$ which contains $t$ variables that are selected at the $t$th iteration. Then $P_t=X_{\Pi_t}(X_{\Pi_t}^HX_{\Pi_t})^{-1}X_{\Pi_t}^H$ is the projection matrix onto the linear subspace spanned by the columns of $X_{\Pi_t}$, and we assume $P_0=0$. 

We want to prove $\Pi_t\subset\Pi$ by induction. First at $t=0$, $\Pi_t=\emptyset\subset\Pi$. Assume at iteration $t$, $\Pi_t\subset\Pi$, then the residual $r_t$ can be written as
$$ r_t = (I-P_t)y = (I-P_t)X_{\Pi}z +(I-P_t)\eta\triangleq s_t+n_t.$$

Let $M_{\Pi}^t = \|X_{\Pi}^H s_t\|_\infty$, $M_{\Pi^c}^ t = \|X_{\Pi^c}^H s_t\|_\infty$ and $N_t =\|X^H n_t \|_{\infty}$, then a sufficient condition for $\Pi_{t+1}\subset\Pi$, i.e. for OMP to select a correct variable at the next iteration is that 
\begin{equation} \label{cond}
M_{\Pi}^t-M_{\Pi^c}^t>2N_t, 
\end{equation}
since
$$  \|X_{\Pi}^H r_t\|_\infty   \geq M_{\Pi}^t -N_t > M_{\Pi^c}^t+N_t \geq  \|X_{\Pi^c}^H r_t\|_\infty.$$

%  which happens with probability at least $1- 2(p^{\alpha}\sqrt{2\pi\log p})^{-1}$ . Define 

\subsection{Bounding $M_{\Pi}^t$ and $M_{\Pi^c}^t$} 
Define the event 
$$\cG_1 = \{X~\mbox{satisfies the}~(k,\epsilon,\delta)\mbox{-StOC} \},$$ 
that happens with probability at least $1-4p^{-1}$ with respect to $\bar{\Pi}$ from Proposition~\ref{StOC}, and the event 
$$\cG_2 =\{  \|X_{\Pi}^HX_{\Pi}-I \|_2\leq 1/2 \},$$ 
which happens at least $1-2p^{-2\log 2}$ with respect to $\bar{\Pi}$ from Proposition~\ref{submat}. Let the event $\mathcal{G}=\cG_1\cap \cG_2$. From the above discussions the event $\mathcal{G}$ holds with probability at least $1-4p^{-1}-2p^{-2\log 2}$ with respect to $\bar{\Pi}$. 

Now we bound $M_{\Pi}^t$ and $M_{\Pi^c}^t$ under the event $\cG$. Let $\Pi_t^c = \Pi\backslash\Pi_t$ be the set of yet to be selected indices of the support of $\beta$, and $\beta_{\Pi_t^c}=z_{\Pi_t^c}$ be the corresponding coefficients. Since $(I-P_t) X_{\Pi}z\in\mathcal{R}(X_{\Pi_t^c})$ belongs to the linear subspace spanned by the columns of $X_{\Pi_t^c}$, we can find a vector $w$ of dimension $(k-t)$ such that $X_{\Pi_t^c}w=(I-P_t) X_{\Pi}z$, where the vector $w$ can be written as 
\begin{align*}
w & =(X_{\Pi_t^c}^HX_{\Pi_t^c})^{-1}X_{\Pi_t^c}^H(I-P_t)X_{\Pi_t^c} z_{\Pi_t^c} \\
&= z_{\Pi_t^c} - (X_{\Pi_t^c}^HX_{\Pi_t^c})^{-1}X_{\Pi_t^c}^HP_tX_{\Pi_t^c} z_{\Pi_t^c} .
\end{align*}

We need the following lemma on eigenvalue majorization from \cite{CaiWang}.
\begin{lem}[\cite{CaiWang}] \label{Cai}
The minimum and maximum eigenvalues of $X_{\Pi_t^c}^H(I-P_t)X_{\Pi_t^c}$ are bounded as
\begin{align*}
\lambda_{\min}(X_{\Pi_t^c}^H(I-P_t)X_{\Pi_t^c}) & \geq \lambda_{\min}(X_{\Pi}^HX_{\Pi}), \\
\lambda_{\max}(X_{\Pi_t^c}^H(I-P_t)X_{\Pi_t^c}) & \leq \lambda_{\max}(X_{\Pi}^HX_{\Pi}).
\end{align*}
\end{lem}
The readers are referred to \cite{CaiWang} for the proof. Since we have
\begin{align}
\|w\|_2 & \leq \|(X_{\Pi_t^c}^HX_{\Pi_t^c})^{-1}\|_2\|X_{\Pi_t^c}^H(I-P_t)X_{\Pi_t^c} z_{\Pi_t^c}\|_2 \label{w1} \\
& \leq 2\|X_{\Pi}^HX_{\Pi} \|_2 \|z_{\Pi_t^c} \|_2 \leq 3 \| z_{\Pi_t^c} \|_2, \label{w2}
\end{align}
where \eqref{w1} follows from Lemma~\ref{Cai}, and \eqref{w2} follows from Proposition~\ref{submat}. Also, 
\begin{align*}
 \|X_{\Pi_t^c}^HP_tX_{\Pi_t^c} z_{\Pi_t^c}\|_{\infty} &= \|X_{\Pi_t^c}^HX_{\Pi_t}(X_{\Pi_t}^HX_{\Pi_t})^{-1}X_{\Pi_t}^HX_{\Pi_t^c} z_{\Pi_t^c}\|_{\infty}  \\
% \leq \|(X_{\Pi_t^c}^HX_{\Pi_t^c})\|_2\|(X_{\Pi_t^c}^HX_{\Pi_t^c})^{-1}X_{\Pi_t^c}^HP_tX_{\Pi_t^c} z_{\Pi_t^c}\|_2 \\
 & \leq \epsilon\|(X_{\Pi_t}^HX_{\Pi_t})^{-1}X_{\Pi_t}^HX_{\Pi_t^c} z_{\Pi_t^c}\|_2 \\
& \leq \epsilon\|(X_{\Pi_t}^HX_{\Pi_t})^{-1} \|_2 \|X_{\Pi_t}^HX_{\Pi_t^c}\|_2 \|z_{\Pi_t^c}\|_2\\
& \leq \epsilon\|z_{\Pi_t^c}\|_2,
  %3/2(\|y\|_2 +  \| z_{\Pi_t^c} \|_2) = 6  \| z_{\Pi_t^c} \|_2.
 \end{align*}
 therefore $M_{\Pi}^t$ can be bounded as
\begin{align}
M_{\Pi}^t &= \| X_{\Pi_t^c}^HX_{\Pi_t^c} z_{\Pi_t^c} - X_{\Pi_t^c}^HP_tX_{\Pi_t^c} z_{\Pi_t^c} \|_\infty \nonumber \\
& \geq \| z_{\Pi_t^c} \|_\infty -\|( X_{\Pi_t^c}^HX_{\Pi_t^c} -I) z_{\Pi_t^c} \|_\infty- \|X_{\Pi_t^c}^HP_tX_{\Pi_t^c} z_{\Pi_t^c} \|_\infty  \nonumber \\
& \geq \|z_{\Pi_t^c} \|_{\infty}- 2\epsilon \| z_{\Pi_t^c} \|_2. \label{M1}
\end{align}
where \eqref{M1} follows from \eqref{StOC-1}. Next, $M_{\Pi^c}^t$ can be bounded as
\begin{align} 
M_{\Pi^c}^ t &= \|X_{\Pi^c}^H (I-P_t) X_{\Pi}z\|_\infty \nonumber \\
& = \|X_{\Pi^c}^H X_{\Pi_t^c} w\|_\infty \nonumber \\
& \leq \epsilon \|w\|_2\leq 3\epsilon \| z_{\Pi_t^c} \|_2. \label{M2}
\end{align}
where \eqref{M2} follows from \eqref{StOC-2}.

\subsection{Bounding $N_t$}
Conditioned on the event $\cG$, for each $P_t$, since $I-P_t$ is also a projection matrix, define the event 
\begin{equation}
\mathcal{H}_t = \{N_t\leq\sigma \sqrt{(1+\alpha)\log p}\},\quad t=0,\cdots, k-1. \label{Nt}
\end{equation}
Then from Proposition~\ref{noise}, it happens with probability at least $1- (p^{\alpha}\pi)^{-1}$ with respect to $\eta$. We further define the event $\mathcal{H}=\cap_{t=0}^{k-1} \mathcal{H}_t$, then from the union bound $\Pr(\cH|\cG) = \Pr(\cH)  \geq 1- k(p^{\alpha}\pi)^{-1}$. Similarly, for the event $\mathcal{H}'=\cap_{t=0}^{k} \mathcal{H}_t$, then from the union bound $\Pr(\cH'|\cG) = \Pr(\cH')  \geq 1- (k+1)(p^{\alpha}\pi)^{-1}$.

\subsection{Proof of Theorem~\ref{thm-omp-k} and Theorem~\ref{thm-omp-stop}}
Define the event $\mathcal{I}=\cG\cap\cH$, from the above discussions we have $\Pr(\mathcal{I})\geq 1-k(p^{\alpha}\pi)^{-1}-2p^{-2\log 2}-4p^{-1}$. Now we are ready to analyze the OMP performance under the event $\mathcal{G}$. We want to prove $\Pi_t\subset\Pi$ by induction.

Substituting the bounds \eqref{M1}, \eqref{M2} and \eqref{Nt} into \eqref{cond}, it is sufficient that at the $t$th iteration
\begin{equation} \label{bound1}
\|z_{\Pi_t^c} \|_{\infty} > 5\epsilon \| z_{\Pi_t^c} \|_2  + 2\sigma\sqrt{ (1+\alpha)\log p}.
\end{equation}

Note that $\|z_{\Pi_t^c}\|_\infty \geq |\beta|_{(t+1)}$, $ \| z_{\Pi_t^c} \|_2 \leq \sqrt{k-t}|\beta|_{(t+1)}$, \eqref{bound1} is satisfied by 
$$  |\beta|_{(t+1)} > 5\epsilon \sqrt{k-t}  |\beta|_{(t+1)} + 2\sigma\sqrt{ (1+\alpha)\log p}. $$
Since $k<1/(c_1^1\mu^2\log p)$ and $\epsilon=10\mu\sqrt{2\log p}$, this is equivalent to the condition in \eqref{lar-condition} for $0\leq t\leq k-1$, therefore a correct variable is selected at the $t$th iteration, $\Pi_t\subset\Pi$. Since the sparsity level of $\beta$ is $k$, the OMP algorithm in Algorithm~\ref{OMP-k} successfully finds the support of $\beta$ in $k$ iterations under the event $\mathcal{I}$, and we have proved Theorem~\ref{thm-omp-k}.

Now we define the event $\cI'=\cH'\cap\cG$, where $\cI'\subset\cI$ which happens with probability at least $1-(k+1)(p^{\alpha}\pi)^{-1}-2p^{-2\log 2}-4p^{-1}$. Conditioned on the event $\cI'$, in order to prove Theorem~\ref{thm-omp-stop}, we need to further show that $\|X^Hr_t\|_\infty> \delta$ for $0\leq t\leq k-1$ so that the algorithm doesn't stop early, and $\|X^Hr_k\|_\infty\leq  \delta$ so that the algorithm stops at the $k$th iteration. While the latter is obvious from the definition of $\cH_k$, for the first inequality we have
\begin{align}
\|X^Hr_t\|_\infty  & \geq M_{\Pi}^t -N_t \nonumber \\
& \geq \|z_{\Pi_t^c} \|_{\infty}- 2\epsilon \| z_{\Pi_t^c} \|_2 -\sigma\sqrt{ (1+\alpha)\log p} \label{stop-1} \\
& > 3\epsilon \| z_{\Pi_t^c} \|_2  +\sigma\sqrt{ (1+\alpha)\log p} \label{stop-2} \\
& \geq \sigma\sqrt{ (1+\alpha)\log p} = \delta, \nonumber
\end{align}
where \eqref{stop-1} follows from \eqref{M1} and \eqref{Nt}, and \eqref{stop-2} follows from \eqref{bound1}. 

\textbf{Remark:} The proof of Corollary~\ref{partial} is straightforward by early-terminating the induction procedure at the $k'$th iteration.

\subsection{Proof of Corollary~\ref{partial-largest}}
Again we prove by induction. First at $t=0$, $\Pi_t=\emptyset\subset\Pi$. Assume at iteration $t$, the OMP algorithm has successfully detected the $t$ largest entries of $|\beta|$. For $i\in\Pi_t^c$ that corresponds to the $t+1$th largest entry of $|\beta|$, we have
\begin{align*}
 |x_i^Hr_t| & \geq   |z_i |- 2\epsilon \| z_{\Pi_t^c} \|_2 -\sigma\sqrt{(1+\alpha)\log p}, \\
 & = |\beta|_{(t+1)}- 2\epsilon \| z_{\Pi_t^c} \|_2 -\sigma\sqrt{(1+\alpha)\log p}
 \end{align*}
 from a simple variation of \eqref{M1}. On the other hand, for $j\in\Pi_t^c$ that corresponds to the rest undetected entries of $\beta$, we have
\begin{align*}
 |x_j^Hr_t| & \leq |z_j| + 2\epsilon \| z_{\Pi_t^c} \|_2 + \sigma\sqrt{(1+\alpha)\log p}, \\
 &  \leq |\beta|_{(t+2)} + 2\epsilon \| z_{\Pi_t^c} \|_2 + \sigma\sqrt{(1+\alpha)\log p}.
\end{align*}
To detect the $t+1$th largest entries it is sufficient to have  
\begin{equation*}
|\beta|_{(t+1)}- |\beta|_{(t+2)} \geq 4\epsilon\sqrt{k-t}|\beta|_{(t+1)} +2\sigma\sqrt{(1+\alpha)\log p}.
\end{equation*}
This is satisfied when \eqref{decay} by simply plugging it into the above equation.

\subsection{Proof of Theorem~\ref{reconstruction}}
Conditioned on the event that the support is successfully recovered, since 
$$\hat{z}=X_{\Pi}^{\dag}y=X_{\Pi}^{\dag}(X_{\Pi}z+\eta)=z+(X_{\Pi}^HX_{\Pi})^{-1}X_{\Pi}^T\eta,$$
where $z=\beta_{\Pi}$ is the non-zero part of $\beta$, and
\begin{align}
&\;\;\;\;\; \| (X_{\Pi}^HX_{\Pi})^{-1} X_{\Pi}^H \eta \|_2^2 \nonumber \\
 &\leq \|(X_{\Pi}^HX_{\Pi})^{-1}\|_2^2 \|X_{\Pi}^H\eta\|_2^2 \nonumber \\
& \leq 4k \|X_{\Pi}^H\eta\|_\infty^2 \leq 4(1+\alpha)k\sigma^2\log p,
\end{align}
it follows that 
$$ \| \hat{\beta}-\beta\|_2^2 = \| \hat{z}-z\|_2^2 \leq 4(1+\alpha)k\sigma^2\log p.  $$

\section{Conclusion} \label{conclusion}
In this paper, we provide coherence-based performance guarantees of two variants of OMP for both support recovery and signal reconstruction of sparse signals when the measurements are corrupted by noise. It is shown that if $X$ satisfied the strong coherence property, then with $n\gtrsim\mathcal{O}(k\log p)$, OMP recovers a $k$-sparse signal with high probability. In particular, the  guarantees obtained here separate the properties required of the measurement matrix from the properties required of the signal. The resilience of OMP to variability in relative strength of the entries of the signal might be an advantage in applications like multi-user detection in wireless communications because it makes power control less critical \cite{XCAC,MUD_journal}.

\section*{Acknowledgements}

The work of Y. Chi and R. Calderbank was supported in part by NSF under Grants NSF CCF-0915299 and NSF CCF-1017431.

\bibliography{omp_ref}
\bibliographystyle{IEEEtran}

\end{document}